\documentclass[12pt,reqno]{article}
\usepackage{amsmath}
\usepackage{amssymb}

\usepackage{latexsym}
\usepackage{amsfonts}
\usepackage{amsmath}
\usepackage{times,pifont}

\newtheorem{theorem}{Theorem}[section]

\newtheorem{proposition}[theorem]{Proposition}
\newtheorem{lemma}[theorem]{Lemma}
\newtheorem{corollary}[theorem]{Corollary}
\newtheorem{definition}[theorem]{Definition}
\newtheorem{conjecture}[theorem]{Conjecture}
\newenvironment{proof}{{\em Proof. }}{\hfill $\Box$ \vspace{1em}}

\newcommand{\E}{{\mathbb E}}
\def\MAJ{\boldsymbol{\mathrm{MAJ}}}

\newcommand{\calT}{{\mathcal T}}

\newcommand{\Var}{{\mathbb V}{\rm ar}}
\newcommand{\Cov}{{\mathbb C}{\rm ov}}
\newcommand{\RR}{{\mathbb R}}

\newcommand{\calC}{{\mathcal C}}

\newcommand{\Z}{{\mathbb Z}}
\newcommand{\T}{{\mathbb T}}
\newcommand{\Pro}{{\mathbb P}}

\def\Inf{{\bf I}}

\newcommand{\one}{{\mathbf 1}}

\newcommand{\bi}{\begin{itemize}}
\newcommand{\ei}{\end{itemize}}
\newcommand{\be}{\begin{enumerate}}
\newcommand{\ee}{\end{enumerate}}
\newcommand{\ra}{\rightarrow}

\newcommand{\ep}{\epsilon}

\newcommand{\iy}{\infty}
\newcommand{\beq}{\begin{equation}}
\newcommand{\eeq}{\end{equation}}
\newcommand{\beqa}{\begin{eqnarray*}}
\newcommand{\eeqa}{\end{eqnarray*}}
\newcommand{\btm}{\begin{theorem}}
\newcommand{\etm}{\end{theorem}}
\newcommand{\bpf}{\begin{proof}}
\newcommand{\epf}{\end{proof}}
\newcommand{\bla}{\begin{lemma}}
\newcommand{\ela}{\end{lemma}}
\newcommand{\bdn}{\begin{definition}}
\newcommand{\edn}{\end{definition}}
\newcommand{\bpn}{\begin{proposition}}
\newcommand{\epn}{\end{proposition}}
\newcommand{\bcy}{\begin{corollary}}
\newcommand{\ecy}{\end{corollary}}

\begin{document}

\title{Volatility of Boolean functions}
\author{Johan Jonasson\thanks{Chalmers University of Technology and University of Gothenburg} \thanks{Research supported by the Knut and Alice Wallenberg Foundation} \thanks{http://www.math.chalmers.se/$\sim$jonasson}
and Jeffrey E. Steif\thanks{Chalmers University of Technology and University of Gothenburg}
\thanks{http://www.math.chalmers.se/$\sim$steif}
\thanks{Research supported by the Knut and Alice Wallenberg Foundation and the Swedish Research Council}}
\maketitle

\begin{abstract}
We study the volatility of the output of a Boolean function when the input bits undergo a natural dynamics.
For $n = 1,2,\ldots$, let $f_n:\{0,1\}^{m_n} \ra \{0,1\}$ be a Boolean function and
$X^{(n)}(t)=(X_1(t),\ldots,X_{m_n}(t))_{t \in [0,\infty)}$ be a vector of i.i.d.\
stationary continuous time Markov chains on $\{0,1\}$ that jump from $0$ to $1$ with rate
$p_n \in [0,1]$ and from $1$ to $0$ with rate $q_n=1-p_n$.  Our object of study will be
$C_n$ which is the number of state changes of $f_n(X^{(n)}(t))$ as a function of $t$ during $[0,1]$.
We say that the family $\{f_n\}_{n\ge 1}$ is volatile if $C_n \ra \iy$ in distribution as
$n\to\infty$ and say that $\{f_n\}_{n\ge 1}$ is tame if $\{C_n\}_{n\ge 1}$ is tight.
We study these concepts in and of themselves as well as investigate their relationship
with the recent notions of noise sensitivity and noise stability. In addition, we
study the question of lameness which means that $\Pro(C_n =0)\ra 1$ as $n\to\infty$.
Finally, we investigate these properties for 
the majority function, iterated 3-majority, the AND/OR function on the binary tree 
and percolation on certain trees in various regimes.
\end{abstract}

\noindent{\em AMS Subject classification :\/ 60K99 } \\
\noindent{\em Key words and phrases:\/  Boolean function, noise sensitivity, noise stability}  \\
\noindent{\em Short title: Volatility of Boolean functions}

\section{Introduction}

We are given a sequence of Boolean functions $\{f_n\}_{n\ge 1}$ with
$f_n:\{0,1\}^{m_n} \ra \{0,1\}$ for some sequence $\{m_n\}$ and also given a sequence $\{p_n\} \in [0,1]$.
Denoting $\{1,2,\ldots, k\}$ by $[k]$, for each $n$ and
for each $i \in [m_n]$, let $\{X^{(n)}_i(t)\}_{t \in [0,\infty)}$ be the stationary continuous time
Markov process on $\{0,1\}$ that jumps from $0$ to $1$ with rate $p_n$ and from $1$ to $0$
with rate $1-p_n$ started in stationarity. (Equivalently $X^{(n)}_i(t)$
updates with rate $1$ and at a given update,
the value is chosen to be 1 with probability $p_n$ and 0 with probability $1-p_n$
independently of everything else.) Assume that the $\{X^{(n)}_i(t)\}_{t\in [0,\infty)}$ are
independent as $i$ varies and write $X^{(n)}(t)$ for $(X^{(n)}_1(t),\ldots,X^{(n)}_{m_n}(t))$.
Finally, the object of our focus will be $C_n([a,b])$ which is defined to be the number of
times that $f_n(X^{(n)}(t))$ changes its state during the time interval $[a,b]$. We abbreviate
$C_n([0,1])$ by $C_n$ throughout.

We say that $\{f_n\}$ is {\em degenerate} with respect to $\{p_n\}$ if
\[
\lim_{n\to\infty}\Pro(f_n(X^{(n)}(0))=1)\Pro(f_n(X^{(n)}(0))=0)=0
\]
and
{\em nondegenerate} with respect to $\{p_n\}$ if for some $\delta >0$,
\[
\Pro(f_n(X^{(n)}(0))=1)\in [\delta, 1-\delta]
\]
for all $n$. (Note that a sequence can of course be neither degenerate nor nondegenerate although it
will always have a subsequence which is either one or the other.)
The first concept we give captures the notion that it is unlikely that there is any change of state.

\begin{definition} \label{DDgen}
We say that $\{f_n\}_{n\ge 1}$ is {\em  lame with respect to $\{p_n\}$}
if
\[
\lim_{n\to\infty}\Pro(C_n)=0)= 1.
\]
\end{definition}

The first relatively easy proposition says that a necessary condition for lameness is that the
sequence is degenerate.

\bpn \label{prop.DDgenImpliesDgen}
Let $\{f_n\}_{n\ge 1}$ be a sequence of Boolean functions and $\{p_n\}$ be a sequence in $[0,1]$.
If $\{f_n\}_{n\ge 1}$ is lame with respect to $\{p_n\}$, then
it is degenerate with respect to $\{p_n\}$.
\epn

The following two definitions will be central to the paper.

\begin{definition} \label{DSe}
We say that $\{f_n\}_{n\ge 1}$ is {\em volatile with respect to $\{p_n\}$}
if $C_n$ approaches $\infty$ in distribution.
\end{definition}

\begin{definition} \label{DSt}
We say that $\{f_n\}_{n\ge 1}$ is {\em tame with respect to $\{p_n\}$}
if $\{C_n\}_{n\ge 1}$ is tight.
\end{definition}

Note that lameness is a special case of tameness.
While it is obvious that the notions of lameness, volatility and tameness
may depend on the sequence $\{p_n\}$, it is natural to ask if these definitions depend on
the length of the time interval chosen which we have taken to be 1. Lameness
clearly does not. The fact that tameness does
not depend on the length of the time interval is straightforward and follows from
the fact that if a sequence of random vectors $(X_n,Y_n)$ is such that $\{X_n\}_{n\ge 1}$ and
$\{Y_n\}_{n\ge 1}$ are each tight, then $\{X_n+Y_n\}_{n\ge 1}$ is also tight.
The fact that
volatility does not depend on the length of the time interval, while certainly believable and in fact
true, does not follow from such general considerations. Rather, some explicit
properties of the process are needed to establish this. It turns out that Markovianness and reversibility are
sufficient. This follows from the following lemma, whose proof presented later
is not so difficult.

\bla \label{lemma.DSwelldefined}
Let $\{f_n\}_{n\ge 1}$ be a sequence of Boolean functions and $\{p_n\}$ be a sequence in $[0,1]$.
The following four conditions are equivalent.\\
(i).  The sequence is volatile.\\
(ii). For each $\delta>0$,  $C_n([0,\delta])$ approaches $\infty$ in distribution. \\
(iii). $\lim_{n\to\infty}\Pro(C_n \ge 1) = 1$.\\
(iv). For all $\delta>0$,  $\lim_{n\to\infty}\Pro(C_n([0,\delta]) \ge 1) = 1$.
\ela

The following corollary, which will be used in a number of our examples, easily follows.

\bcy \label{la}
Assume that $\lim_{n\to\infty}\Pro(f_n(X^{(n)}(0))=1) = 0$ and that
\[
\lim_{n\to\infty}\Pro(\exists t \in [0,1]:f_n(X^{(n)}(t))=1) =1.
\]
Then $\{f_n\}_{n\ge 1}$ is volatile.
\ecy

It is trivial to construct a sequence of functions which is neither volatile nor tame;
simply let $\{p_n\}\equiv 1/2$, $f_n$ be the dictator function (which simply outputs the value of the
first bit) for even $n$ and the parity function (which outputs the mod 2 sum of the values of
all the bits) for odd $n$. The following definition captures the notion of a sequence being neither
volatile nor tame but for less trivial reasons.

\begin{definition} \label{DSmixed}
We say that $\{f_n\}_{n\ge 1}$ is {\em semi-volatile with respect to $\{p_n\}$}
if
\[
\liminf_n \Pro(C_n=0) >0 \mbox{ and } \lim_{M\to\infty}\liminf_{n\to\infty}\Pro(C_n>M) >0.
\]
\end{definition}

This means that for all large $n$, the distribution of $C_n$ has some weight at 0 and some weight near
$\infty$. It is elementary to check that $\{f_n\}_{n\ge 1}$ is semi-volatile with respect to $\{p_n\}$
if and only $\{f_n\}_{n\ge 1}$ does not have any subsequence which is either volatile
with respect to $\{p_n\}$ or tame with respect to $\{p_n\}$.

A simple example of a sequence of Boolean functions on $n$ bits which is semi-volatile when
$\{p_n\}\equiv 1/2$ is the function which is 1 if both the value of the first bit is 1
and the mod 2 sum of the values of the remaining bits is 0.
(This function is known as dictator AND parity.)

It turns out to be natural, with hindsight, to further partition the class of 
semi-volatile sequences into two groups.
For lack of a better name, we call them Type 1 and Type 2.

\begin{definition} \label{DSmixedTYpe1}
We say that $\{f_n\}_{n\ge 1}$ is {\em Type 1 semi-volatile with respect to $\{p_n\}$}
if it is semi-volatile and if for all $k\ge 1$
\[
\lim_{n\to\infty}\Pro(1\le C_n\le k) =0.
\]
\end{definition}

\begin{definition} \label{DSmixedTYpe2}
We say that $\{f_n\}_{n\ge 1}$ is {\em Type 2 semi-volatile with respect to $\{p_n\}$}
if it is semi-volatile and if for some $k\ge 1$
\[
\limsup_{n\to\infty}\Pro(1\le C_n\le k) >0.
\]
\end{definition}

It is an elementary exercise to check that the ``dictator AND parity'' example
given above is Type 1 semi-volatile.
Type 1 behavior can be viewed as a discrete time analogue of,
for example, the fact that if we start a Brownian motion from
1, the number of times during $[0,1]$ that it crosses 0 is either 0 or infinite.
While we believe that most ``easy'' examples of semi-volatile functions are Type 1, it is not hard to construct an example having Type 2.

\bpn \label{prop.Type2}
There exists a nondegenerate sequence of Boolean functions which is Type 2 semi-volatile.
\epn

\noindent
{\bf Remark.} We will see that later results also yield 
Type 2 semi-volatility in some given situations.

\smallskip

The above concepts are related to, but distinct from, the recent notions of noise sensitivity and
noise stability. These two latter concepts were introduced in the seminal paper by
Benjamini, Kalai and Schramm (\cite{BKS}) and were developed further in \cite{GS}.
We now give the relevant definitions.

In the following, we are given a sequence of integers $\{m_n\}$ and numbers
$p_n\in [0,1]$ and consider the product measure on
$\{0,1\}^{m_n}$ with marginal distribution $p_n\delta_1+(1-p_n)\delta_0$.
We denote a random element of $\{0,1\}^{m_n}$ under this measure by $\omega_n$ and given $\epsilon >0$, we let
$\omega_n^\epsilon$ denote the configuration obtained from $\omega_n$
where each bit of $\omega_n$ is independently with probability $\epsilon$ replaced by a 1 or 0,
with respective probabilities $p_n$ and $1-p_n$ independently of everything else.

\begin{definition} \label{NS}
The sequence of Boolean functions $f_{n}:\{0,1\}^{m_n}\rightarrow \{0,1\}$
is said to be \emph{noise sensitive} w.r.t.\ $\{p_n\}$ if for any $\epsilon>0$,
\begin{equation} \label{eqn:NS}
\lim_{n\rightarrow \infty}
\E[f_n(\omega_n)f_n(\omega_n^\epsilon)]- \E[f_n(\omega_n)]^2 =0.
\end{equation}
\end{definition}

\begin{definition} \label{defstab1}
The sequence of Boolean functions $f_{n}:\{0,1\}^{m_n}\rightarrow \{0,1\}$
is said to be noise stable w.r.t.\ $\{p_n\}$
if for any $\delta>0$ there exists an $\epsilon>0$ such that
\[
\sup_n \Pro(f_{n}(\omega_n)\neq f_{n}(\omega_n^{\epsilon}))\le \delta\,.
\]
\end{definition}

\bigskip\noindent
A trivial but key observation relating noise sensitivity/stability to the concepts introduced earlier
is that
\begin{equation} \label{eqn:SameJoint}
(X^{(n)}(0),X^{(n)}(t)) \mbox{ and } (\omega_n,\omega_n^{1-e^{-t}})
\end{equation}
have the same joint distribution.

We now state our result relating the notions of noise stability and tameness.

\bpn \label{prop.DStabImpliesStab}
Let $\{f_n\}_{n\ge 1}$ be a sequence of Boolean functions and $\{p_n\}$ be a sequence in $[0,1]$.
If $\{f_n\}_{n\ge 1}$ is tame with respect to $\{p_n\}$, then it is noise stable with respect to $\{p_n\}$.
\epn

The majority function, defined next, yields a simple example which is noise stable but not tame.

\medskip
\begin{definition}[Majority function]\label{def.maj}
Let $n$ be odd and define
\[
\MAJ_n(x_1,\ldots,x_n):= I_{\{\sum_{i=1}^n x_i \ge (n+1)/2\}}.
\]
\end{definition}

\bpn \label{prop.MajorityNotDStab}
With $\{p_n\}\equiv 1/2$ and only considering odd $n$, the sequence $\{\MAJ_n\}_{n\ge 1}$ is noise stable
but not tame.
\epn

\bigskip
As will be remarked later after the proof of this result, this sequence is in fact Type 1 semi-volatile.
A much more striking example is the following model, which yields very rich behavior.

The model we want to consider is dynamical percolation; see \cite{S} for a survey of this subject.
We shall be brief here and we will appeal to results in
\cite{HPS} and \cite{PSS} as well as to results in \cite{Ly} which concern ordinary percolation.
We consider an infinite tree which is spherically symmetric which
means that all vertices at a given level have the same number of children (which may depend on
the level). We will perform percolation on the edges of our tree with parameter $1/2$ which means
each edge is independently removed with probability $1/2$. An edge is considered in state 1 if it is
retained and 0 if it is removed. We now let $f_n$ be the Boolean function of these variables which
is 1 if there is a path from the root to the $n$th level using only the retained edges and 0
otherwise.

\btm \label{theorem.StableButDsens}
Consider the above case of dynamical percolation and let $w_n$ denote the expected number of
paths from the root to the $n$th level which only use retained edges.
For the sequence of associated Boolean functions $\{f_n\}_{n\ge 1}$ defined above,
we have the following. \\
(i). For any sequence $w_n$, $\{f_n\}_{n\ge 1}$ is noise stable. \\
(ii). If $w_n \asymp \log n$, then
$\lim_{n\to\infty}\Pro(\forall t\in [0,1]:\,\,f_n(X^{(n)}(t))=0)=1$ and hence
$\{f_n\}$ is lame. (Degeneracy follows from by Proposition \ref{prop.DDgenImpliesDgen}.)\\
(iii). If $w_n \asymp (\log n)^{1+\delta}$ for some $\delta >0$, then
$\lim_{n\to\infty}\Pro(f_n(X^{(n)}(0))=1)=0$ (implying that $\{f_n\}$ is degenerate) and
$\{f_n\}$ is Type 1 semi-volatile. \\
(iv). If $w_n \asymp n(\log n)^\alpha$ for $\alpha \in (1,2]$, then $\{f_n\}$ is nondegenerate and
Type 1 semi-volatile. \\
%Moreover,
%$\lim_{n\to\infty}\Pro(\forall t\in [0,1]:\,\,f_n(X^{(n)}(t))=0)>0$
%while $\lim_{n\to\infty}\Pro(\forall t\in [0,1]:\,\,f_n(X^{(n)}(t))=1)=0$.
(v).  If $w_n \asymp n(\log n)^\alpha$ for $\alpha>2$, then $\{f_n\}$ is nondegenerate and Type 2
semi-volatile.\\
%Moreover $\lim_{n\to\infty}\Pro(\forall t\in [0,1]:\,\,f_n(X^{(n)}(t))=0)>0$
%and $\lim_{n\to\infty}\Pro(\forall t\in [0,1]:\,\,f_n(X^{(n)}(t))=1)>0$.
(vi).  If $w_n \asymp n^\alpha$ for $\alpha>2$, then $\{f_n\}$ is nondegenerate and tame.
%Moreover $\lim_{n\to\infty}\Pro(\forall t\in [0,1]:\,\,f_n(X^{(n)}(t))=0)>0$
%and $\lim_{n\to\infty}\Pro(\forall t\in [0,1]:\,\,f_n(X^{(n)}(t))=1)>0$.
\etm

\medskip\noindent
{\bf Remark.} For each of the cases, one can construct a tree so that 
$w_n$ fulfills the stated assumption.
The case (iv) is very different from the cases (v) and (vi) as they correspond to very different behaviors
in the dynamical percolation models. Namely, in the regime of (iv), there exist exceptional times at which
there are no infinite clusters while in the regimes of (v) and (vi), there are no such exceptional times.
This difference turns out to correspond to
$\lim_{n\to\infty}\Pro(\forall t\in [0,1]:\,\,f_n(X^{(n)}(t))=1)$ being 0 in case (iv) and positive in cases
(v) and (vi).

The next result gives us our main relationship between noise sensitivity and volatility.

\bpn \label{NSImpliesNSense}
Let $\{f_n\}_{n\ge 1}$ be a sequence of Boolean functions and $\{p_n\}$ be a sequence in $[0,1]$.
If $\{f_n\}_{n\ge 1}$ is nondegenerate and noise sensitive with respect to $\{p_n\}$, then it is volatile with respect to $\{p_n\}$.
\epn

\medskip\noindent
{\bf Remark.} The nondegeneracy condition is needed since a degenerate
sequence is immediately noise sensitive and certainly might not be volatile,
for example if $f_n$ is the constant function 1 for each $n$.

\medskip
Without using this language, the implication in Proposition \ref{NSImpliesNSense}
was proved in \cite{BKS} for the specific sequence of Boolean
functions corresponding to percolation crossings of a large square in $\Z^2$.
Their proof technique however proves the above more
general result. This is presented in \cite{GS} (again without using this language) when $p_n=1/2$.
However, the proof in \cite{GS} does not in fact require this latter assumption on $p_n$. We therefore
give no proof of Proposition \ref{NSImpliesNSense} but rather refer the reader to either Corollary 5.1 in
\cite{BKS} or Chapter 1 of \cite{GS}.

Returning to tameness, there is a nice sufficient condition in terms of influences.
The notion of influence, which we now introduce, is crucial in the study of noise sensitivity and
noise stability. To explain this notion, we first endow $\{0,1\}^{n}$
with product measure $\Pro_p$ with marginal distribution
$p\delta_1+(1-p)\delta_0$ and consider a given Boolean function $f$ defined on
$\{0,1\}^{n}$. We again denote a typical element of $\{0,1\}^{n}$ by $\omega_n$ and given $i\in [n]$,
we let $\omega_n^i$ be $\omega_n$ but with $i$ rerandomized to be 1 or 0, with respective
probabilities $p$ and $1-p$.

\begin{definition} \label{influence}
The {\em influence} of bit $i$ on $f$ at parameter $p$, denoted by $\Inf^p_i(f)$, is
$\Pro_p(f(\omega_n)\neq f(\omega_n^i))$. The quantity $\sum_i\Inf^p_i(f_n)$ is referred to as
the {\em total influence} at parameter $p$.
\end{definition}

\bpn \label{prop.FiniteInfluenceYieldsD-Stab}
Let $\{f_n\}_{n\ge 1}$ be a sequence of Boolean functions and $\{p_n\}$ be a sequence in $[0,1]$.
Then for each $n$, $\E[C_n]=\sum_i \Inf^{p_n}_i(f_n)$. Consequently,
if $\sup_n \sum_i \Inf^{p_n}_i(f_n)<\infty$, then $\{f_n\}_{n\ge 1}$ is tame with respect to $\{p_n\}$.
\epn

The majority functions, $\{\MAJ_n\}_{n\ge 1}$, $p_n \equiv 1/2$, show that $\sup_n \sum_i \Inf^{p_n}_i(f_n)<\infty$, while known to be
sufficient, is not a necessary condition for noise stability since it is easy to check that
$\sum_i \Inf^{p_n}_i(f_n)$ is of order $\sqrt{n}$ in this case while noise stability is well known.
One might ask if $\sup_n \sum_i \Inf^{p_n}_i(f_n)<\infty$ is however necessary for the stronger property of tameness;
i.e., whether the converse of the last statement of 
Proposition \ref{prop.FiniteInfluenceYieldsD-Stab} might be true. The answer
turns out to be no as stated next.

\bpn \label{prop.BigInfluenceButDStab}
There exists a sequence of Boolean functions $\{f_n\}_{n\ge 1}$
which is nondegenerate with respect to $\{p_n\}\equiv 1/2$,
satisfies $\lim_{n\to\infty}\E[C_n]=\infty$ but is tame.
\epn

It is reasonable to expect that the converse of
Proposition \ref{prop.FiniteInfluenceYieldsD-Stab} perhaps holds under some reasonable
additional conditions. We point out that using standard second moment methods and the
Paley-Zygmund inequality, it is standard to check that
$\lim_{n\to\infty}\E[C_n]=\infty$ and
$\E[C^2_n]\le O(1)\E[C_n]^2$ implies that the sequence is not
tame while $\lim_{n\to\infty}\E[C_n]=\infty$ and
$\E[C^2_n]\le (1+o(1))\E[C_n]^2$ implies that the sequence is volatile.
Therefore, one approach to establishing volatility or non-tameness
for various classes of functions would be to show that
$\lim_{n\to\infty}\E[C_n]=\infty$ and attempt to obtain good estimates
on the second moment of $C_n$.

Concerning the possibility that the converse of the last statement of 
Proposition \ref{prop.FiniteInfluenceYieldsD-Stab}
holds under some reasonable assumptions, we have the following conjecture.

\begin{conjecture}
Let $\{f_n\}_{n\ge 1}$ be a sequence of Boolean functions which is nondegenerate with 
respect to $\{p_n\}$
satisfying $\lim_{n\to\infty}\E[C_n]=\infty$. %% and $\{p_n\}$ bounded away from $0$ and $1$. 
If each $f_n$ is transitive (meaning that there is
a transitive group action on $[m_n]$ which leaves $f_n$ invariant), then
the sequence is not tame with respect to $\{p_n\}$.
\end{conjecture}

\medskip\noindent
{\bf Remarks.} (i). If there exists $\delta>0$ such that
$\delta\le p_n\le 1-\delta$ for each $n$, then, using the main result in \cite{KKL},
transitivity and nondegeneracy implies that $\lim_{n\to\infty}\E[C_n]=\infty$.  
(Technically, \cite{KKL} only covers the case when $p=1/2$ and 
\cite{BKKKL} is needed for general $p$. \\
%%(ii). The above implication is not true without 
%%the assumption on the sequence $\{p_n\}$ as illustrated
%%by letting $f_n$ be the event that there is at least one 1 and $p_n=1/n$. \\
(ii). The conjecture is false if one drops the nondegeneracy assumption 
but of course keeping the
$\lim_{n\to\infty}\E[C_n]=\infty$ assumption. This can be seen by letting $p_n\equiv 1/2$ and
let $f_n$ be the event that there are at least $n/2 +\sqrt{n}c_n$ 1's where $c_n$ 
increases very slowly to infinity. In this case, one in fact has lameness.

\medskip

The following important result, due to Benjamini, Kalai and Schramm \cite{BKS}, relates the sum of the squared influences to noise sensitivity.

\btm \label{thm.squared_influences}
If $\{p_n\}$ is bounded away from $0$ and $1$ and  $\sum_i\Inf_i^{p_n}(f_n)^2 \ra 0$ as $n \ra \iy$, then $\{f_n\}$ is noise sensitive.
\etm

\bigskip
We now discuss a situation where the sequence $\{f_n\}_{n\ge 1}$ is monotone in a certain sense. (This is note the usual definition of a Boolean function being monotone.)

\begin{definition} \label{monotone}
The sequence $\{f_n\}_{n\ge 1}$ is {\em monotone} if for any $\omega\in \{0,1\}^{m_n}$,
$f_n(\omega)=1$ implies that $f_{n-1}(\omega')=1$ where $\omega'$ is $\omega$ restricted to the first
$m_{n-1}$ bits.
\end{definition}

A key example of a monotone sequence is the sequence of Boolean functions
treated in Theorem \ref{theorem.StableButDsens}. Note that while all previous definitions
are unaffected if one changes 0 and 1, the above definition is affected.
The following proposition due to Erik Broman, while not hard, is of interest to point out.

\bpn \label{prop.Broman} (E. Broman)
Any monotone sequence of Boolean functions is noise stable.
\epn

When we have a monotone sequence of Boolean functions $\{f_n\}_{n\ge 1}$, we obtain in a natural way
a function $f_\infty$ on the space $\{0,1\}^\infty$ defined by
\[
f_\infty(\omega):=\lim_{n\to\infty} f_n(\omega_n)
\]
where $\omega_n$ is $\omega$ restricted to the first $m_n$ bits. (Monotonicity of course implies
the existence of the limit.) Since $X^{(\infty)}(t)=(X_1(t),\ldots)_{t \in [0,\infty)}$ has obvious
meaning, we can consider the process
\begin{equation} \label{fInfinity}
\{f_\infty(X^{(\infty)}(t))\}_{t \in [0,1]}.
\end{equation}
In this situation, the behavior of the various dynamical properties that we have been
studying can be expressed in terms of the process given in  (\ref{fInfinity}).
For example, consider the degenerate situation when
$\lim_{n\to\infty}\Pro(f_n(X^{(n)}(0))=1)=0$ which is equivalent to $\Pro(f_n(X^{(\infty)}(0))=1)=0$.
It can then be shown for example that lameness is equivalent (assuming we are not
in the trivial case $f_n\equiv 1$ for all $n$) to
$\Pro(\forall t\in [0,1]\,\,f_\infty(X^{(\infty)}(t))=0)=1$.
The latter property, in slightly different language, was studied in \cite{BHPS}. Namely, one has a
measurable function $f_\infty$ from $\{0,1\}^\infty$ into $\{0,1\}$ (not necessarily a limit of
functions as above) with $\Pro(f_\infty(X^{(\infty)}(0))=0)=1$ and one asks if it is also the case that
$\Pro(\forall t\in [0,1]\,\,f_\infty(X^{(\infty)}(t))=0)=1$. When the opposite
\[
\Pro(f_\infty(X^{(\infty)}(t))=1\,\,  \mbox{ for some } t\in [0,1])>0
\]
holds, we say $f_\infty$ is called {\em dynamically sensitive} and refer to the times $t$ for which
$f_\infty(X^{(\infty)}(t))=1$ as {\em exceptional times}.
The question of dynamical sensitivity was posed and answered in \cite{BHPS} for a number of functions
$f_\infty$ which corresponded to various standard concepts in probability theory (such as strong law of large
numbers, a.s.\ central limit theorems, recurrent/transience, run lengths, etc.). In the general Markov process
lingo, a property which is dynamically sensitive is often called {\em nonpolar}.

We end this introduction by looking at two particular (sequence of) Boolean functions that we will
analyze it in some detail.  They are the Iterated 3-majority function and the 
AND/OR function on the binary tree.

We first deal with the Iterated 3-majority function.
Let $\T=\T_n$ be the rooted ternary tree of depth $n$. To each leaf, $\l$, attach an independent Bernoulli random variable $Y_l$, the {\em state} of $l$, with $\Pro(Y_l=1)=p=p_n$.
The states of the other vertices are recursively defined by
\[Y_u=\max(Y_{v_1}Y_{v_2},Y_{v_1}Y_{v_3},Y_{v_2}Y_{v_3}),\]
where $v_1,v_2,v_3$ are the children of $u$,
i.e.\ $Y_u$ is defined to be the majority of its children.
We are interested in $f(Y)=Y_o$, where $o$ is the root; this is the so called {\em iterated 3-majority} function.

If $p=1/2$, then obviously $\Pro(f(Y)=1)=1/2$. It is also known (see \cite{AS}) that if $p=1/2+\gamma(2/3)^n$ for $\gamma = \Theta(1)$, then $\Pro(f(Y)=1)=(1/2)(1+r)$ where $r$ is bounded away from $0$ and $\pm 1$, with $r \ra \pm 1$ whenever $\gamma \ra \pm \iy$.
By recursion one can readily see that for any $p$, the influence, $I_l$, of a leaf variable, $Y_l$, is at most $2^{-n}$. Hence $\sum_l I_l^2 \ra 0$ and it follows from Theorem \ref{thm.squared_influences} that $f$ is noise sensitive.

Dynamics is introduced to the model in the usual way:
let the leaf variables $Y_l$ update according to unit intensity Poisson processes as usual. So, in the general setup, we have that the $X^{(n)}_i(t)$ are the dynamical leaf variables, $i=1,\ldots,3^n$.
For $p=1/2-\ep$, with $\ep=\gamma(2/3)^n$ and $\gamma=\Theta(1)$, we just noted that $f$ is noise sensitive. As $f$ is non-degenerate (see \cite{AS}), it follows from
Proposition \ref{NSImpliesNSense} also volatile.
The natural question here is: if $\gamma = \gamma_n \ra \iy$, then for what $\gamma$ will $f$ still be volatile, or, equivalently by Corollary \ref{la}, for which $\gamma$ will there be exceptional times at which $f(X(t))=1$ even though for fixed $t$,
$\Pro(f(X(t))=0) \ra 1$?
The following result shows that $\gamma$ polynomial in $n$ (or logarithmic in the number of variables if you like) is the interesting range of orders and that there is a sharp cutoff.

\btm \label{tf}
In the setup above, let $\gamma = \gamma_n = n^\alpha$ for a constant $\alpha$ and $p=p_n=1/2-\gamma(2/3)^n$. If $\alpha>\alpha_0:=\log(3/2)/\log(2)$, then $f$ is lame and if $\alpha<\alpha_0$, then $f$ is volatile.
\etm

Lastly we deal with the AND/OR function on the binary tree.
Let $\T_n$ be the rooted binary tree of depth $n$. Regard this as an electric network where at each vertex, there is either an AND-gate or an OR-gate.
Then supplying two 0/1 in-signals at each leaf gives a certain out-signal from the root.
We will assume that each gate is chosen to be AND or OR independently with probability $1/2$ and each leaf gets one 0 in-signal and one 1 in-signal.
(Equivalently we could have let the in-signals also be random, but in this way all the randomness goes into the choice of the states of the gates.)
By the symmetry $f(x)=1-f(1-x)$, the out-signal at the root is 1 with probability $1/2$ and $0$ with probability $1/2$.
Now introduce dynamics as above on the states of the gates. Let $X_v(t)$ be the process of states for the gate at vertex $v$ (leaves not included), $X(t)=X^{(n)}(t)=\{X^{(n)}_v(t)\}$ and let $f(X(t))=f_n(X^{(n)}(t))$ be the out-signal at the root.
We prove the following.

\btm \label{te}
The out-signal at the root for the dynamical AND/OR-process on the binary tree is Type 2 semi-volatile with respect to $p \equiv 1/2$.
In particular, using symmetry, $\Pro(\forall t:f(X(t))=1)$ is bounded away from $0$ and $1$.
\etm

\section{Proofs of general results}

\medskip\noindent
{\em Proof of Proposition  \ref{prop.DDgenImpliesDgen}.} \,\,
It is easily seen that it suffices to show that for any $\delta>0$, there exists a $\epsilon>0$ so that
for any $n$, $p_n$ and Boolean function $f$ on $\{0,1\}^{n}$ satisfying $\Pro(f=1)\in [\delta, 1-\delta]$,
we have that
$\Pro(f(X^{(n)}(0))\neq f(X^{(n)}(1)))\ge\epsilon$.
The display following (2.3) in \cite{LS} and Equation \eqref{eqn:SameJoint} here together yield that
\[
\Cov(f_n(X^{(n)}(0)),f_n(X^{(n)}(1)))\le (1-e^{-1})\Var(f_n(X^{(n)}(0))).
\]
The assumption that $\Pro(f=1)\in [\delta, 1-\delta]$ implies that
\[
\Var(f_n(X^{(n)}(0)))\ge \delta(1-\delta).
\]
This and the previous display easily yield that
\[
\Pro(f(X^{(n)}(0))\neq f(X^{(n)}(1)))\ge\epsilon
\]
for some $\epsilon$, only depending on $\delta$.
\hfill $\Box$

\medskip\noindent
{\em Proof of Lemma \ref{lemma.DSwelldefined}.} \,\,
We will show (i) implies (ii) and (iii) implies (iv). This suffices since (ii) implies (iii)
is vacuous while (iv) implies (i) is elementary and left to the reader.

\noindent
{\bf (i) implies (ii).} 
For this, it suffices to show that for any $\delta >0$, $C_n([0,2\delta])$
approaches $\infty$ in distribution implies that $C_n([0,\delta])$ approaches $\infty$ in distribution
which we now argue. If the latter is not true, then there exists $M$ and $\epsilon_0$ so that for
infinitely many $n$
\[
\Pro(C_n([0,\delta])\le M)\ge \epsilon_0.
\]
This yields that for infinitely many $n$,
\[
\E[\Pro(C_n([0,\delta])\le M\mid X^{(n)}(\delta))]\ge \epsilon_0.
\]
It is easy to see that if $g$ is a function on a probability space taking values in
$[0,1]$ with $\int g \ge \epsilon_0$, then $\Pro(g\ge \epsilon_0/2)\ge \epsilon_0/2$.
Hence, for infinitely many $n$, there is a subset $T_n\subseteq \{0,1\}^{m_n}$ so that
$\Pro(T_n)\ge \epsilon_0/2$ and
\[
\min_{\eta\in T_n} \Pro(C_n([0,\delta])\le M\mid X^{(n)}(\delta)=\eta)\ge \epsilon_0/2.
\]
By Markovianness and time reversibility, we have that for such $n$, for $\eta \in T_n$,
\[
\Pro(\{C_n([0,\delta])\le M\}\cap \{C_n([\delta,2\delta])\le M\}\mid X^{(n)}(\delta)=\eta)\ge \epsilon^2_0/4
\]
and hence
\[
\Pro(C_n([0,2\delta])\le 2M\mid X^{(n)}(\delta)=\eta)\ge \epsilon^2_0/4.
\]
Since $\Pro(T_n)\ge \epsilon_0/2$, we obtain that for infinitely many $n$
\[
\Pro(C_n([0,2\delta])\le 2M)\ge \epsilon^3_0/8,
\]
implying that $C_n(2\delta)$ does not approach $\infty$ in distribution.

\noindent
{\bf (iii) implies (iv).} 
It suffices, by iteration, to show that for any $a$, \\
$\lim_{n\to\infty}\Pro(C_n([0,a]) \ge 1) = 1$ implies that
$\lim_{n\to\infty}\Pro(C_n([0,\frac{a}{2}]) \ge 1) = 1$. The former implies that
\[
\lim_{n\to\infty}\E[\Pro(C_n([0,a])\ge 1\mid X^{(n)}(\frac{a}{2}))]= 1.
\]
It follows that there exist subsets $T_n\subseteq \{0,1\}^{m_n}$ so that $\Pro(T_n)\ra 1$ and
\[
\lim_{n\to\infty}\inf_{\eta\in T_n}\Pro(C_n([0,a])\ge 1\mid X^{(n)}(\frac{a}{2})=\eta)]= 1.
\]
If $A_n,B_n$ are independent events with $\Pro(A_n)=\Pro(B_n)$ and
$\Pro(A_n\cup B_n)\ra 1$, it follows that $\Pro(A_n\cap B_n)\ra 1$.
Applying this with $A_n:=\{C_n([0,\frac{a}{2}])\ge 1\}$ and $B_n:=\{C_n([\frac{a}{2},a])\ge 1\}$, yields that
\[
\lim_{n\to\infty}\inf_{\eta\in T_n}\Pro(C_n([0,\frac{a}{2}])\ge 1\mid X^{(n)}(\frac{a}{2})=\eta)]= 1.
\]
Since $\lim_{n\to\infty}\Pro(T_n)= 1$, we obtain $\lim_{n\to\infty}\Pro(C_n([0,\frac{a}{2}]) \ge 1) = 1$,
as desired.
\hfill $\Box$

\medskip\noindent
{\em Proof of Proposition \ref{prop.Type2}.} \,\,
Let $f_n$ be the function which outputs the value of the first bit if the value of the second bit is 1
and which outputs the mod 2 sum of the values of bits $3,\ldots,n$ if the value of the second bit is 0.
It is elementary to check that this example has the desired properties.
\hfill $\Box$

\medskip\noindent
{\em Proof of Proposition \ref{prop.DStabImpliesStab}.} \,\,
Let $\delta >0$. Choose an integer $M$ so that
\[
\sup_n\Pro(C_n\ge M)< \delta/2.
\]
Let $k\ge \frac{2M}{\delta}$ be an integer. It is easily seen by symmetry %and Markov's inequality
that for any $n$
\[
\Pro(C_n([0,\frac{1}{k}])>0\mid  C_n\le M)\le \frac{M}{k}.
\]
This easily yields from the above that
\[
\sup_n \Pro(C_n([0,\frac{1}{k}])>0)\le \frac{\delta}{2} +\frac{M}{k}\le \delta
\]
and hence
\[
\sup_n \Pro(f_n(X^{(n)}(0))\neq f_n(X^{(n)}({\frac{1}{k}})))\le \delta.
\]
By (\ref{eqn:SameJoint}), this easily yields noise stability.
\hfill $\Box$

\medskip\noindent
{\em Proof of Proposition \ref{prop.MajorityNotDStab}.} \,\,
The fact that this sequence is noise stable is relatively standard; see \cite{BKS} or \cite{GS}.
It is well known that
\[
\large\{\frac{\sum_{k=1}^n2X^{(n)}_k(t)-1}{\sqrt{n}}\large\}_{t\ge 0}
\]
converges in distribution (with respect to the Skorohod topology) to the Ornstein-Uhlenbeck process,
which can be described as the stationary Gaussian process
$\{U(t)\}_{t\ge 0}$ with continuous paths having mean zero and the convariance structure
$\E(U(t) U(t+s))=e^{-s}$. It is known that on finite time intervals, $\{U(t)\}_{t\ge 0}$ and Brownian
motion are absolutely continuous with respect to each other and from this it is easy to
see that $\{U_t\}$ crosses 0 infinitely many times in $[0,1]$ with positive probability. By the above
convergence, it is easy to show that
\[
\lim_{M\to\infty}\liminf_{n\to\infty}\Pro(C_n>M) >0.
\]
This rules out tameness.
\hfill $\Box$

\medskip\noindent
{\bf Remark.} Since with positive probability $\{U_t\}$ never crosses 0 during $[0,1]$, it also follows that
$\inf_n \Pro(C_n=0) >0$ which also rules out volatility. Therefore this sequence
is in fact semi-volatile and it is not hard to show that it is Type 1.

\medskip\noindent
{\em Proof of Theorem \ref{theorem.StableButDsens}.} \,\,
(i). This sequence of Boolean functions is clearly monotone in the sense of
Definition \ref{monotone} and hence the noise stability follows from
Proposition \ref{prop.Broman} (to be proved later).  \\
(ii). We first mention that while degeneracy follows from the first statement,
degeneracy is a direct consequence of Theorem 2.1 in \cite{Ly} which implies that
$\lim_{n\to\infty}\Pro(f_n(X^{(n)}(0))=1)=0$. The first statement is
a consequence of Theorem 1.5 in \cite{HPS} as well as its proof.

We postpone moving to the remaining cases since the following discussion is relevant to all of these cases.
We let
$$
\calT_1:=\{t\in [0,1]: f_n(X^{(n)}(t))=1\,\,\forall n\},
$$
$$
\calT_0:=[0,1]\backslash \calT_1= \{t\in [0,1]: f_n(X^{(n)}(t))=0 \mbox{ for some n}\},
$$
and
$$
\calT^n_1:=\{t\in [0,1]: f_n(X^{(n)}(t))=1\}.
$$

We first claim that $\calT_1$ is either empty or infinite (in fact uncountable). This follows from
the following general theorem explained to us by Steve Evans. For any stationary
reversible Markov process, the set of times in $[0,1]$ at which the process is in a certain subset of the state space
is either empty or infinite. The proof of this result is detailed in Lemma 2.3 in \cite{PS}.
In our particular case of percolation and where $\{f_n\}$ is degenerate, a hands on proof
that $\calT_1$ is either empty or infinite is given in Lemma 3.4 in \cite{HPS}.

We next let $\calC$ denote the number of components of $\calT_1$ in the interval $[0,1]$. It is elementary to check that
if $\calC\ge 2$, then
\[
\liminf_{n\to\infty} C_n \ge \calC.
\]
It follows from Fatou's lemma that for $M\ge 2$,
\[
\liminf_{n\to\infty} \Pro(C_n\ge M)\ge  \Pro(\calC\ge M)
\]
and hence that
\begin{equation} \label{eqn:CompBigGivesFlipsBig}
\lim_{M\to\infty}\liminf_{n\to\infty} \Pro(C_n\ge M)\ge  \Pro(\calC=\infty).
\end{equation}

By the earlier result above, up to a set of measure 0, our probability space is partitioned into
\begin{equation} \label{eqn:Partition}
\{\calT_1= \emptyset\}\cup\{\calT_1= \infty,\calC =\infty\}\cup\{\calT_1= \infty,\calC <\infty\}.
\end{equation}

We observe that in all cases,
$\inf_n\Pro(\forall t\in [0,1]:\,\,f_n(X^{(n)}(t))=0)>0$ since with positive
probability the edges emanating from the root are off throughout $[0,1]$.  This yields that
$\liminf_n \Pro(C_n=0) >0$.

Therefore, if the middle event in Equation \eqref{eqn:Partition} has positive probability, then semi-volatility
follows from Equation \eqref{eqn:CompBigGivesFlipsBig}.  We next argue that if
the middle event in Equation \eqref{eqn:Partition} has positive probability and the
last event in Equation \eqref{eqn:Partition} has 0 probability, then we have Type 1 semi-volatility.

To see this, we first claim that
\begin{equation} \label{eqn:FirstCase}
\liminf_n \Pro(C_n=0)\ge \Pro(\calT_1= \emptyset).
\end{equation}
To see this, we have $\Pro(C_n=0)\ge \Pro(\calT^n_1= \emptyset)$ yielding
$\liminf_n \Pro(C_n=0)\ge \lim_n \Pro(\calT^n_1= \emptyset)$. Finally, while it is not true
that for every $\omega$ $\{\calT_1= \emptyset\}=\cup_n\{\calT^n_1= \emptyset\}$, it is true that
these two latter events are the same
up to a set of measure 0 by Lemma 3.2 in \cite{HPS}. This yields Equation \eqref{eqn:FirstCase}.

We finally observe that if the first two events in Equation \eqref{eqn:Partition} have positive probability but
the third has probability 0, then, by Equations \eqref{eqn:CompBigGivesFlipsBig} and  \eqref{eqn:FirstCase},
we have that
\[
\liminf_n \Pro(C_n=0)+\lim_{M\to\infty}\liminf_{n\to\infty} \Pro(C_n\ge M)=1.
\]
It is elementary to see that this implies that once one has established semi-volatility,
then Type 1 semi-volatility follows. We now return to the remaining cases.

(iii). Theorem 2.1 in \cite{Ly} again implies that
$\lim_{n\to\infty}\Pro(f_n(X^{(n)}(0))=1)=0$ and hence the sequence is degenerate.
This also implies that a.s.\ $\calT_0$ is dense which in turn implies that the third event in
Equation \eqref{eqn:Partition} has probability 0.
Next, a consequence of Theorem 1.5 in \cite{HPS} is that $\Pro(\calT_1 \neq\emptyset)>0$
which implies that the second event in Equation \eqref{eqn:Partition} has positive probability.
From the above discussion, we can conclude Type 1 semi-volatility.\\
(iv). Theorem 2.1 in \cite{Ly} also implies in this case that
$\lim_{n\to\infty}\Pro(f_n(X^{(n)}(0))=1)=0$ and hence the sequence is degenerate.
The nondegeneracy of $\{f_n\}$ (together with the monotonicity of the sequence) implies,
by Fubini's Theorem, that the expected value of the Lebesgue measure of $\calT_1$ is positive and
hence with positive probability, $\calT_1$ is infinite.
Next Theorem 1.2(ii) in \cite{PSS} together with Kolmogorov's 0-1 Law implies that
a.s.\ $\calT_0$ is dense. Exactly as in case (iii), this now proves Type 1 semi-volatility.\\
(v). Theorem 2.1 in \cite{Ly} implies again in this case that $\{f_n\}$ is nondegenerate.
The fact that
\begin{equation} \label{eqn:GoodAllTheWay}
\lim_{n\to\infty}\Pro(\forall t\in [0,1]:\,\,f_n(X^{(n)}(t))=1)>0
\end{equation}
follows from
Theorem 1.1(i) and Lemma 3.2 in \cite{PSS}. Theorem 1.3(ii) in \cite{PSS} says that $\calC$ is infinite
with positive probability. From the discussion before (iii), this implies semi-volatility.
We now show the sequence is Type 2 semi-volatile.
Let $A$ be the event that the first child of the root has at all
times in $[0,1]$ an infinite open path to infinity (not going through the root) and that the edges
to the other children are always off during $[0,1]$. Then $A$ has positive probability by
Equation \eqref{eqn:GoodAllTheWay} and, in this case, the number of state changes is
precisely the number of times that the edge from
the root to the first child switches. This yields Type 2. \\
(vi).  Theorem 2.1 in \cite{Ly} implies again in this case that $\{f_n\}$ is nondegenerate.
Theorem 1.3(ii) in \cite{PSS} says that $C<\infty$ a.s.
This suggests tameness but does not immediately imply it. However the proof of Theorem 1.3(ii) in \cite{PSS}
shows that $\E[C_n]=O(1)$ yielding tameness.
\hfill $\Box$

\medskip\noindent
{\em Proof of Proposition \ref{prop.FiniteInfluenceYieldsD-Stab}.} \,\,
Letting $C_n(i)$ be the number of times that there is a state change due to a flip at location $i$,
it suffices to show that for each
$i$, $\Inf^{p_n}_i(f_n) =\E(C_n(i))$. The number of times that location $i$ is rerandomized has a Poisson
distribution with mean 1. The probability that a particular rerandomization at location $i$ causes a
state change is precisely $\Inf^{p_n}_i(f_n)$. This easily yields the equality. The final statement holds since
a uniform bound on the expectation of a set of random variables immediately yields tightness.
\hfill $\Box$

\medskip\noindent
{\em Proof of Proposition \ref{prop.BigInfluenceButDStab}.} \,\,
Let $m_n:= 1+n+3^n$ and $f_n$ be the Boolean function which is
the first bit unless bits 2 through n+1 are all 1's, in which case $f_n$ outputs the mod 2 sum
of the last $3^n$ bits. Clearly the sequence is nondegenerate
and it is elementary to check that the total influence is at least $(3/2)^n$ but
that the sequence is tame.
\hfill $\Box$

\medskip\noindent
{\em Proof of Proposition \ref{prop.Broman}.} \,\,
Let $\delta>0$. Since $f_n(\omega)$ converges for every $\omega$, we can choose $N$ so that for all $n\ge N$,
$\Pro(f_n(\omega)\neq f_N(\omega))<\delta/4$. Next, choose $\epsilon >0$ so that
$\Pro(f_N(\omega)\neq f_N(\omega^\epsilon))<\delta/2$. Then, for all $n\ge N$,
$\Pro(f_n(\omega)\neq f_n(\omega^\epsilon))$ is at most
\[
\Pro(f_n(\omega)\neq f_N(\omega))+
\Pro(f_N(\omega)\neq f_N(\omega^\epsilon))
+\Pro(f_N(\omega^\epsilon)\neq f_n(\omega^\epsilon)) < \delta.
\]
By choosing $\epsilon$ even smaller we can insure that
$\Pro(f_k(\omega)\neq f_k(\omega^\epsilon)) < \delta$ for $k< n$ as well. This completes the proof.
\hfill $\Box$
\bigskip

\section{Proof for iterated 3-majority}

{Proof of Theorem \ref{tf}.} \,\,
Let us start by remarking that here, as elsewhere, most of the numbers appearing are in no way optimal.
Fix $n$, $\gamma = \gamma_n =n^\alpha$ and $\ep = \ep_n = \gamma (2/3)^n$.
With $\alpha$ fixed, everything from here on is true for $n$ sufficiently large.
Let $a_k$ be the probability that the root of a $k$-generation ternary tree is in state $1$. Then $\{a_k\}$ satisfies the recursion
\beq \label{eba}
a_{k+1} = 3a_k^2 - 2a_k^3,\, a_0=p
\eeq
Note that $n=\log_{3/2}(\gamma/\ep)$.
We want to compute a good approximation of $a_n$.
The idea is to use that for (roughly) the first $\log_{3/2}(1/\ep)$ steps of the recursion, $a_k$ stays close to $1/2$, but from this point it drops rapidly for the remaining (roughly) $\log_{2/3}(\gamma)$ steps.
Writing $a_k=(1/2)(1-\pi_k)$, (\ref{eba}) gives
\beq \label{ebb}
\pi_{k+1}=\frac32 \pi_k-\frac12 \pi_k^3,\, \pi_0=2\ep
\eeq
As long as $\pi_k$ is small, we have that ``$\pi_k$ is very close to $2\ep(3/2)^k$''.
To quantify this, let $\pi'_0=\pi_0$ and $\pi'_{k+1} = (3/2)\pi'_k$.
Then $\pi'_k = 2\ep(3/2)^k$. Obviously $\pi_k < \pi'_k$ for $k \geq 1$ and hence $\pi_k < 1/2$ for $k \leq \log_{3/2}(1/4\ep)$.
Writing $\pi'_k = c_k \pi_k$, we have, since $\pi_k$ is increasing in $k$ by (\ref{ebb}),
\beqa
c_{k+1}\pi_{k+1} &=& c_k \pi_{k+1}+\frac12 c_k \pi_k^3 \\
&=& c_k\left(1+\frac12\frac{\pi_k^3}{\pi_{k+1}}\right) \pi_{k+1} \\
&<& c_k(1+\frac12\pi_k^2)\pi_{k+1},
\eeqa
so $c_{k+1} < (1+(1/2)\pi_k^2)c_k$.
Hence, as long as $\pi_k<1/2$, $c_{k+1}/c_k<9/8$ and hence $\pi_{k+1}/\pi_k > (3/2)(8/9) = 4/3$, i.e.\ as long as $\pi_k <1/2$, $\pi_k$ increases in $k$ at least exponentially at rate $4/3$.
From this it follows that, as long as $\pi_k<1/2$,
\[c_k < \prod_{j=0}^{k-1}\left(1+\frac12\pi_j^2\right) < \prod_{l=0}^{\iy}\left(1+\frac12\left((1/2)(3/4)^{l}\right)^2\right) < \frac43.\]

It now follows that $\pi_k>1/2$ for $k \geq \log_{3/2}(1/4\ep)+1$, since
if $\pi'_j \geq 1/2$ and $\pi_j<1/2$, then $\pi_j>3/8$, since $c_j<4/3$,
so that $\pi_{j+1} > (3/2)(3/8)-(1/2)(3/8)^3 > 1/2$.

However as soon as $a_k$ starts to deviate significantly from $1/2$, it also starts to decay very rapidly.
To quantify this, let $d_{k+1}=3d_k^2$, $d_0=1/4$.
This can be solved explicitly:
\[d_k = \frac13\left(\frac34\right)^{2^k} = \frac13\left(\frac23\right)^{2\xi 2^k},\]
where $2\xi = \log(4/3)/\log(3/2)$.
Now obviously $a_{j+k} \leq d_k$, $k \geq 0$, i.e.\ $a_k \leq d_{k-j}$, $k \geq j$, for any $j$ such that $a_j \leq 1/4$.
If also $a_j > 5/32$, then $a_k \geq d_{k-j+2}$.
To see this, write $e_{k+1}=3e_k^2-2e_k^3$, $e_0=5/32$ and $e_k=g_k d_k$.
Then we get
\[g_{k+1} = g_k^2-\frac23 g_k^3d_k \geq \frac56 g_k^2,\, g_0=\frac58,\]
where the inequality follows from that $d_k \leq 1/4$.
This gives
\[g_k \geq \left(\frac65\right)\left(\frac{25}{48}\right)^{2^k} > \left(\frac65\right)\left(\frac12\right)^{2^k}.\]
Hence
\[e_k \geq \left(\frac25\right)\left(\frac38\right)^{2^k} > d_{k+2},\]
since $2/5>1/3$ and $(3/4)^4<3/8$.

By the recursion for $a_k$, there is exactly one $j$ for which $5/32 < a_j \leq 1/4$.
We saw above that $j:=\lceil \log_{3/2}(1/4\ep) \rceil +1$ or $j-1$ is such a number.
Since $n-j$ or $n-j+1$ then equals $\lfloor \log_{3/2}(4\gamma) \rfloor$, it follows that
\beqa
a_n &=& \frac13\left(\frac23\right)^{\eta 2^{\log_{3/2}(4\gamma)}} \\
&=& \frac13\left(\frac23\right)^{\eta (4\gamma)^{\alpha_0^{-1}}},
\eeqa
where $\eta=\eta_n \in [\xi/2,8\xi]$.
If $\alpha>\alpha_0$, then the exponent of the right hand side is $n^\nu$ for $\nu>1$ and the right hand side is hence of smaller order than $(1/3)^n$. Since the expected number of updates during unit time is $3^n$, a first moment bound proves the lameness part.

Now we consider the volatility part.
With $n$ again fixed, let $b_k$ be the probability that the root of a depth $k$ tree is in state 1 at times $0$ and $t$, where $t=t_n$ is assumed to be of larger order than $\ep=\gamma(2/3)^n$ (this  assumption will not enter the considerations until a bit further down).
Then $b_0=(1+e^{-t})/4-\ep+\ep^2(1-e^{-t})$.
Now, a depth $k+1$ tree is in state 1 at both times if and only if at least two
of the three depth $k$ subtrees of the root are in state 1 at both times or exactly one of the three subtrees is in state 1 at both times and one of the others is in state 0 at time 0 and state 1 at time $t$ and the other one is vice versa.
This gives the recursion
\beq \label{ebc}
b_{k+1} = 3b_k^2-2b_k^3+6b_k(a_k-b_k)^2.
\eeq
Writing $b_k = a_k(1-\beta_k)$, it follows from (\ref{ebc}) that
\beqa
(1-\beta_{k+1})a_{k+1} &=& 3a_k^2(1-2\beta_k+\beta_k^2) - 2a_k^3(1-3\beta_k+3\beta_k^2-\beta_k^3) \\
& & + 6a_k^3\beta_k^2(1-\beta_k) \\
&=& a_{k+1}-6\beta_k(a_k^2-a_k^3)+3a_k^2\beta_k^2 - 4a_k^3\beta_k^3 \\
&=& a_{k+1}(1-2\beta_k)+2\beta_k a_k^3 + 3\beta_k^2a_k^2 - 4a_k^3\beta_k^3.
\eeqa
It follows that
\[\beta_{k+1} = (2-\frac{2a_k^3}{a_{k+1}})\beta_k - \frac{3a_k^2}{a_{k+1}}\beta_k^2 + \frac{4a_k^3}{a_{k+1}}\beta_k^3,\, \beta_0=(1/2+\ep)(1-e^{-t}).\]
Since $a_k^3/a_{k+1} = p_k/(3-2a_k) \leq 1/4$, (because $a_k \leq 1/2$) and similarly $a_k^2/a_{k+1} \leq 1/2$, we have
\[\beta_{k+1} \geq \frac32\beta_k-\frac32\beta_k^2.\]
Let us compare $\beta'_k$ with $\beta''_k$, where $\beta'_0=\beta''_0=\beta_0$, $\beta'_{k+1}=(3/2)\beta'_k-(3/2)(\beta'_k)^2$ and $\beta''_{k+1}=(3/2)\beta''_k$, in the same manner as for $\pi_k$ and $\pi'_k$ above.
Write $\beta''_k = h_k \beta'_k$. Then
\[h_{k+1}\beta'_{k+1} = \frac32 h_k\beta'_k = h_k\beta'_{k+1}+\frac32 h_k(\beta'_k)^2 \leq \left(1+\frac32 \beta'_k\right)h_k \beta'_{k+1}\]
so $h_{k+1} \leq (1+3\beta'_k/2)h_k$.
Now, as long as $\beta'_k \leq 1/9$, we then have that $\beta'_k$ grows in $k$ at least exponentially at rate $4/3$, which entails that
\[h_k \leq \prod_{l=0}^\iy\left(1+(1/6)(3/4)^l\right) < 2.\]
From this it follows that if $\beta''_k \geq 1/4$, say, then $\beta'_k > 1/9$.
However if $\beta'_k>1/9$, then $\beta'_{k+4} >1/4$.
Hence $\beta'_k$, and hence $\beta_k$, exceeds $1/4$ no later than four steps after $\beta''_k$ exceeds $1/4$. Since the latter happens after $\lceil \log_{3/2}(\beta_0^{-1}/4) \rceil$ steps, it follows in particular that
$\beta_k \geq 1/4$ whenever $k \geq k_0 := \lfloor \log_{3/2}\beta_0^{-1} \rfloor + 2$.
Since $\beta_0$ is of larger order than $\ep$ (since $t$ is of larger order than $\ep$), $\pi_k$ is after $k_0$ steps still very small,
i.e.\ $a_k$ is very close to $1/2$ and will remain so until step
$k_1:=\lceil \log_{3/2}(1/\ep) \rceil -16$, whereas $b_k$ has started to decay rapidly by step $k_0$.
At this point $\pi_{k_1} \leq \pi'_{k_1} < 1/100$, so that $a_{k_1} > 49/100$ and hence $a_{k_1}/a_{k_1+1}<50/49$.
To quantify the rapid decay of $b_k$ from step $k_0$ to $k_1$,
write $b_k=a_k^2(1+\rho_k)$, $k \geq k_0$, insert in (\ref{ebc}) and use (\ref{eba}) to get
\beqa
(1+\rho_{k+1})a_{k+1}^2 &=& 3a_k^4(1+\rho_k)^2 - 2a_k^6(1+\rho_k)^3 + 6a_k^4(1+\rho)(1-a_k-a_k\rho_k)^2 \\
&=& 9a_k^4-12a_k^5+4a_k^6 + (12a_k^4-24a_k^5+12a_k^6)\rho_k \\
& & + (3a_k^4-12a_k^5+12a_k^6)\rho_k^2 + 4a_k^6\rho_k^3 \\
&=& a_{k+1}^2 + 12a_k^4(1-a_k)^2\rho_k + 3a_k^4(1-4a_k+4a_k^2)\rho_k^2 + 4a_k^6 \rho_k^3.
\eeqa
It follows that
\begin{equation} \label{exa}
\rho_{k+1} = \frac{12a_k^4(1-a_k)^2}{a_{k+1}^2}\rho_k + \frac{3a_k^4(1-2a_k)^2}{a_{k+1}^2}\rho_k^2+\frac{4a_k^6}{a_{k+1}^2}\rho_k^3.
\end{equation}
Note that since $\beta_{k_0} \geq 1/4$, we have $a_{k_0}^2(1+\rho_{k_0}) = b_{k_0} \leq (3/4)a_{k_0} < 3/8$,
so $\rho_{k_0} < 3/(8a_{k_0}^2 0.49^2) < 9/16$.
Using this, that $a_k > 49/100$ and $a_k/a_{k+1} < 50/49$ for $k \leq k_1$, it follows from (\ref{exa}) that for $k \leq k_1$,
\[\rho_{k+1} < \frac78 \rho_k.\]
This means that
\beqa
\rho_{k_1} &<& \left(\frac78\right)^{k_1-k_0} \\
&\leq& \left(\frac87\right)^{18}\left(\frac78\right)^{\log_{3/2}(\beta_0/\ep)} \\
&=& \left(\frac87\right)^{18} \left(\frac{\ep}{\beta_0}\right)^c
\eeqa
where $c=\log(8/7)/\log(3/2)$.
For the final steps of the recursions, when $a_k$ is also decaying rapidly, one can extract from (\ref{exa}) and (\ref{eba}) that in any case, as long as $\rho_k<1$
\[\rho_{k+1} < \frac43\rho_k+\frac13\rho_k^2+\frac14\rho_k^3 < 2\rho_k.\]
so for the remaining $\log_{3/2}(\gamma)+16$ steps, $\rho_k$ increases compared to $\rho_{k_1}$ by a factor at most $2^{16}\gamma^d$ where $d=\log(2)/\log(3/2)$.
Multiplying this with the right hand side of the bound for $\rho_{k_1}$, it follows that
\[\rho_n \leq \left(\frac87\right)^{18}2^{16}\gamma^d \left(\frac{\ep}{\beta_0}\right)^c.\]
Recall that $\beta_0 = (1/2+\ep)(1-e^{-t}) = (1+o(1))t/2$. Let $\delta>0$ be an arbitrarily small constant and $t \geq \delta a_n = \delta(1/3)(2/3)^{\eta(4\gamma)^{\alpha_0^{-1}}}$.
We then get that
$\beta_0$ is of larger order than $\ep$ if $\gamma=n^\alpha$ with $\alpha<\alpha_0$ (recall that we needed this assumption earlier).
Moreover, since $\ep=\gamma(2/3)^n$ and
\[\beta_0 > \frac13 t \geq \frac{\delta}{9}\left(\frac23\right)^{\eta(4\gamma)^{\alpha_0^{-1}}} = \frac{\delta}{9}\left(\frac23\right)^{\mu n^{\alpha \alpha_0^{-1}}}\]
with $\mu :=4^{\alpha_0^{-1}}\eta$, we get
\[\frac{\ep}{\beta_0} < \frac{9\gamma}{\delta} \frac{(2/3)^n}{(2/3)^{\mu n^{\alpha\alpha_0^{-1}}}} = \frac{9n^\alpha}{\delta} \frac{(2/3)^n}{(2/3)^{\mu n^{\alpha\alpha_0^{-1}}}}.\]
Hence
\[\rho_n < 3^{18}n^{\alpha d}\left(\frac{9n^\alpha}{\delta} \frac{(2/3)^n}{(2/3)^{c\mu n^{\alpha\alpha_0^{-1}}}}\right)^c < \delta\]
for sufficiently large $n$,
so that
\[b_n = (1+\delta)a_n^2.\]
Now let $Z_n$ be the number of times of the form $j\delta p_n$, $j=0,\ldots, (\delta p_n)^{-1}$ that $f(X(t))=1$.
Then $\E[Z_n]=\delta^{-1}$ and we have just shown that for sufficiently large $n$, $\E[Z_n^2] < \E[Z]+\delta\E[Z_n]^2 < (1+\delta)\E[Z_n]^2$.  It follows that for for sufficiently
large $n$, $\Pro(Z_n>0) > (1+\delta)^{-1}$.
Since $\delta$ was arbitrary, this proves that $\Pro(Z_n>0) \ra 1$ and hence that $f$ is volatile by Corollary \ref{la}.
\hfill $\Box$

\section{Proof for the AND/OR-process on the binary tree}
{\em Proof of Theorem \ref{te}.}\,\,
We start by proving that (an arbitrary subsequence of) $f$ is not tame.
Write $T_l$ and $T_r$ for the left and right subtrees of the root of $\T_n$, respectively (i.e.\ the subtrees that have the left and right children of the root as their roots, respectively).
Let $a^n_k$, $n=0,1,2,\ldots$, $k=0,1,2,\ldots,n$ be the probability that, at a fixed time, $f=1$ and that a given vertex of $\T_n$ at level $k \leq n$ is in state OR and that if that vertex were changed to an AND, that would change $f$ to $0$.
In short, $a^n_k$ is the probability that $f=1$ and the vertex under consideration is pivotal.
Considering the root as generation $0$, we get $a^0_0=1/2$ and $a^n_0=1/4$ for $n \geq 1$ and recursively $a^n_{k+1}=a^{n-1}_k/2$, which can be argued as follows; assuming that our vertex at level $k+1$ is in $T_l$, it is then pivotal if it is pivotal in $T_l$ and either the root is OR and the $T_r$ has out-signal 0 or the root is AND and $T_r$ has out-signal 1.
From this recursion, it follows that $a^n_k=1/2^{k+2}$ for $k=0,\ldots,n-1$ and $a^n_n=1/2^{n+1}$.
This means that if a vertex of $\T_n$ is chosen uniformly at random and updated, the probability that this causes $f$ to go from $1$ to $0$ is $(1/(2N_n))(\sum_{k=0}^{n-1} 2^k(1/2^{k+2}) + 1/2)= (n+2)/(8N_n)$, where $N_n=2^{n+1}-1$ is the number of vertices of $\T_n$.

Let $\Lambda_n$ be the point process consisting of times in $[0,1]$ where $f_n$ changes from $1$ to $0$ and let $S_n$ be 
the total number of points of $\Lambda$.
Recall that the one dimensional intensity function for a general point process $\Lambda$ is the function $\rho_1:[0,1] \ra \RR_+$ such that
$\Pro(\Lambda \cap [x,x+\ep] \neq \emptyset) = \rho_1(x)\ep + o(\ep)$ for all $x$.
The two dimensional intensity function is the function $\rho_2:[0,1] \times [0,1] \ra \RR_+$ such that\
$\Pro(\Lambda \cap [x,x+\ep] \neq \emptyset, \Lambda \cap [y,y+\ep] \neq \emptyset) = \rho_2(x,y)\ep^2 + o(\ep^2)$.
One has in general that $\E[|\Lambda|] = \int_0^1 \rho_1(x)dx$ and $\E[|\Lambda|^2] = \E[|\Lambda|] + \int_0^1\int_0^1 \rho_2(x,y)dx\,dy$.
By what we have shown above together with the fact the the updating intensity of our process is $N_n$,  we get for our point process that $\rho_1(x) = (n+2)/8$.
It follows that $\E[S_n] = (n+2)/8$.

Next we bound the second moment of $S_n$ by bounding $\rho_2$.
For that, we first need to estimate the probability $x_n(t):=\Pro(f_n(X^{(n)}(0))=f_n(X^{(n)}(t))=1)$, $n \geq 0$, for a given $t \in [0,1]$.
Letting $2\tau=1-e^{-t}$, so that $t/4 \leq \tau \leq t/2$, and conditioning on the state of the root at times $0$ and $t$, we get the recursion $x_0(t)=(1-\tau)/2$ and
\beqa
x_{n+1}(t) &=& \frac12(1-\tau)\left(2x_n(t)-x_n(t)^2+2\left(\frac12-x_n(t)\right)^2\right) + \frac12(1-\tau)x_n(t)^2 \\
& & +\tau\left(x_n(t)^2+2x_n(t)\left(\frac12-x_n(t)\right)\right) \\
&=& (1-\tau)(x_n(t)^2+\frac14) + \tau\left(x_n(t)-x_n(t)^2\right).
\eeqa
This recursion can easily be shown to have a unique attractive fixed point at $x=(1/2)(1-\sqrt{\tau}+O(t))$, which is thus approximately $x_n(t)$ is for large $n$.
However, we need an upper bound for $x_n(t)$.
Writing $x_n(t) = (1/2)(1-y_n(t))$, the recursion translates to
\[y_{n+1}(t) = y_n(t)+\frac12 \tau - \frac12 y_n(t)^2 - \tau y_n(t) + \tau y_n(t)^2 := h_\tau(y_n(t)),\, y_0(t)=\tau.\]
Since $\tau<1/2$, we have $h'_\tau(y_n(t)) \in (0,1)$ for all $n$. Since also $h_\tau(\tau)=h_\tau(y_0(t))>y_0(t)$ for any $\tau \in (0,1/2)$, we have that $y_0(t)$ is to the left of the unique positive fixed point of the recursion and it follows that $y_n(t)$ is increasing in $n$ and all $y_n(t)$ are to the left of the fixed point.
As long as $y_n(t) \leq \sqrt{\tau}/12$, we have $y_{n+1}(t) \geq y_n(t) + \tau/2-\tau/144 - \tau^{3/2}/12 > y_n(t) + \tau/3$.
Hence for all $n$ such that $y_n(t) \leq \sqrt{\tau}/12$, $y_n(t) \geq n\tau/3$.
Now if it were the case that $y_n(t) \leq \sqrt{\tau}/12$ for some $n \geq 1/(2\sqrt{\tau})$, we would have also $y_n(t) \geq \sqrt{\tau}/6$, a contradiction.
Hence $y_n(t) > \sqrt{\tau}/12 > \sqrt{t/4}/12 = \sqrt{t}/24$ whenever $n \geq 1/(2\sqrt{\tau})$, which holds whenever
$t \geq 1/n^2$.
Consequently, letting $\beta_n(t) = \one_{(0,1/n^2)}(t)+(1-\sqrt{t}/24)\one_{[1/n^2,1)}(t)$, we have:
\beq \label{er}
x_n(t) \leq \frac12\beta_n(t).
\eeq
(Note that the particular numbers here and later are by no means optimal.)
Now condition on that the system updates at time $0$ and at time $t$ (i.e.\ at times $0$ and $t$ there is some vertex that updates). Let $b_n(t)$ be the conditional probability that $f$ switches from 1 to 0 at both occasions.
We are now going to derive a recursive inequality for $b_n(t)$. Obviously $b_n(t) \leq 1$ for $n = 0,1,2$. Now assume $n \geq 2$ and consider $b_{n+1}(t)$.
One can have a change from $1$ to $0$ at both occasions in the following four distinct ways depending on what two vertices were updated.
\bi
\item[(i)] the root was updated at both times,
\item[(ii)] the root was updated once but not twice,
\item[(iii)] the two updated vertices were taken both from $T_l$ or both from $T_r$,
\item[(iv)] the two updated vertices were taken one from $T_l$ and one from $T_r$.
\ei
Bounding the probability for each of these four ways, using (\ref{er}), and summing gives
\begin{equation} \label{ebn}
b_{n+1}(t) \leq  \frac{1}{16N_{n+1}^2}+\frac{1}{16}\frac{n+2}{N_n N_{n+1}}+\frac{1}{4}\beta_n(t)b_n(t)+\frac{1}{128}\frac{(n+2)^2}{N_n^2}.
\end{equation}
Here the first term comes from the event that the root is chosen for updating at both occasions; the probability for picking the root twice is $(1/N_{n+1})^2$ and given this, the probability that the system changes from $1$ to $0$ at the first occasion is $1/8$ and in addition to this it is required that the root becomes AND after the second update.
The second term comes from (ii); the probability that the root was chosen once but not twice is less than $2/N_{n+1}$. Given such a choice, say that the root was updated at time $0$
and one vertex of $T_l$ was updated at time $t$, for the system to change from $1$ to $0$ on both occasions, it is required that the root goes from OR to AND at time $0$, which has conditional
probability $1/4$, and that $T_l$ goes from $1$ to $0$ at time $t$, which is conditionally independent and has conditional probability $(n+2)/(8N_n)$.
The third term arises from (iii). To see this, note first that the probability for picking the two vertices from the same subtree is less than $1/2$. Then for the whole tree to change from  1 to 0 at both times, of course requires that
the subtree of the chosen vertices also makes the same changes and also that the other subtree has out-signal 0 and the root is OR or
the other subtree has out-signal 1 and the root is AND, at both times. Now this latter thing happens at time 0 with probability $1/2$ and given this, then it also happens at time $t$ if either the root
and the out-signal of the other subtree are the same as at time 0 or if they have both changed. The conditional probability for this
is $(1-\tau)2x_n(t) + \tau(1-2x_n(t))$. Since $\tau<1/2$, this expression is increasing in $x_n(t)$ and is hence no larger than
$(1-\tau)\beta_n(t)+\tau(1-\beta_n(t)) \leq \beta_n(t)$.
Summing up, the contribution from (iii) is at most $\beta_n(t)b_n(t)/4$.
The last term arises from (iv); the probability of choosing vertices from different subtrees is less than $1/2$ and given this the two subtrees must independently have a change from $1$
to $0$ at times $0$ and $t$.
This finishes the proof of (\ref{ebn}).

Bounding the constants in (\ref{ebn}) generously and using that $n+2 \leq 2n$, gives
\[b_{n+1}(t) < \frac14 \beta_n(t)b_n(t) + \frac{n^2}{4^{n+1}}.\]
Using induction it follows that for $n \geq 3$
\beqa
b_n(t) &<& \frac{n^2}{4^{n-1}}(1+\beta_{n-1}(t)+\beta_{n-1}(t)\beta_{n-2}(t)+\ldots+\beta_{n-1}(t)\beta_{n-2}(t)\cdots\beta_2(t)) \\
&<& 32\left(\frac{n}{N_n}\right)^2\sum_{k=2}^{n}\prod_{i=k}^{n-1} \beta_i(t) \\
&<& 32\left(\frac{n}{N_n}\right)^2\sum_{j=0}^\iy \left(1-\frac{\sqrt{t}}{24}\right)^j + \frac{1}{\sqrt{t}}\left(1-\frac{\sqrt{t}}{24}\right)^{(n-1/\sqrt{t})_+} \\
&<& 800\left(\frac{n}{N_n}\right)^2\frac{1}{\sqrt{t}}.
\eeqa
Now the relation between $b_n(t)$ and $f_2$ is, since the process updates with intensity $N_n$, that
$\rho_2(x,y) = N_n^2 b_n(|x-y|) = O(1)(n^2/\sqrt{|x-y|})$.
Integrating over $[0,1] \times [0,1]$, we get 
\[\E[S_n^2] = \E[S_n] + O(1)n^2 = O(1)n^2.\]

Hence $\E[S_n^2]=O(1)\E[S_n]^2$ and the Paley-Zygmund inequality together with the fact that $\E[S_n] \ra \iy$, implies that $f$ cannot be tame.

To show that (an arbitrary subsequence of) $\{f_n\}$ is not volatile, we will now prove that $\Pro(\forall t:f(X(t))=1)$ is bounded away from $0$.
Again we use recursion via conditioning on the root.
Let $T=T_n$ be the first time that the out-signal at the root is $0$ and let $G_n(t) = \Pro(T_n>t)$.
Now, if the root is not OR throughout, $T_n>x$ will occur if the out-signals from the two subtrees of the root both have out-signal $1$ throughout $[0,x]$.
Also, if the root is OR throughout $[0,x]$, then $T_n>x$ occurs if either one of the two subtrees signals $1$ throughout, or if there is some $t \in (0,x)$ such that
one of the subtrees signals $1$ throughout $[0,t]$ but not $[0,x]$ and the other signals $1$ throughout $[t,x]$ but not $[0,x]$.
This gives the recursive inequality
\beqa
G_{n+1}(x) &\geq& \left(1-\frac12e^{-x/2}\right)G_n(x)^2 \\
&+& \frac12e^{-x/2}\left(2G_n(x)-G_n(x)^2 + \int_0^x(G_n(x-t)-G_n(x))(-G_n'(t))dt\right). \\
&\geq& \frac12\left(1+\frac12 x\right)G_n(x)^2 \\
&+& \frac12\left(1-\frac12 x\right)\left(2G_n(x)-G_n(x)^2 + \int_0^x(G_n(x-t)-G_n(x))(-G_n'(t))dt\right).
\eeqa
In more probabilistic notation this becomes
\beqa
G_{n+1}(x) &\geq& \frac12\left(1+\frac12 x\right) \Pro(X \geq x, X' \geq x) \\
&+& \frac12\left(1-\frac12 x\right)\Pro(X+X' \geq x)
\eeqa
where $X$ and $X'$ are independent random variables distributed according to $G_n$. From this, we see that the right hand side is increasing in $G_n$ (since increasing $G_n$ corresponds to making $X$ and $X'$ stochastically larger).
However inserting $(1/2)(1-4\sqrt{x})$ for $G_n(x)$ on the right hand side, one finds that the expression becomes at least as large as $(1/2)(1-4\sqrt{x})$ for all $x \in (0,1)$.
Since $G_0(x)>(1/2)(1-4\sqrt{x})$, this proves that $G_n(x) > (1/2)(1-4\sqrt{x})$ for all $n$ and $x \in (0,1)$ and that no subsequence of $\{f_n\}$ is volatile. Since no subsequence of $\{f_n\}$ was tame, semi-volatility is established by definition.
It remains to prove Type 2. This follows on observing that each of the following events occur independently and with probability bounded away from $0$:
\bi
\item[(i)] the out-signal at the left child of the root is $0$ throughout,
\item[(ii)] the out-signal at the right child of the root is $1$ throughout,
\item[(iii)] the state of the root changes $k$ times.
\ei
\hfill $\Box$

\section*{Acknowledgements}
We thank Erik Broman for allowing us to include his Proposition \ref{prop.Broman} and Steve Evans for some
discussions.

% \noindent
% Johan Jonasson \\
% Department of Mathematics \\
% Chalmers University of Technology and G\"oteborg University \\
% 412 96 \, Gothenburg, Sweden

% \medskip \noindent
% Jeffrey E. Steif \\
% Department of Mathematics \\
% Chalmers University of Technology and G\"oteborg University \\
% 412 96 \, Gothenburg, Sweden
\end{document}